\documentclass[10pt,a4paper,fleqn]{article}
\usepackage[utf8]{inputenc}
\usepackage{amsmath}
\usepackage{amsfonts}
\usepackage{amssymb}
\usepackage{xcolor}
\usepackage{graphicx}
\usepackage{subcaption}
\usepackage{float}
\numberwithin{equation}{section}

\author{Giulio Morpurgo}
\title{On a stricter Twin Primes Conjecture, and on the Polignac's Conjecture in general}
\begin{document}
\maketitle
\begin{abstract}
The Polignac's Conjecture, first formulated by Alphonse de Polignac in 1849, asserts that, for any even number $M$, there exist infinitely many couples of prime numbers $P$, $P+M$. When $M = 2$, this reduces to the Twin Primes Conjecture. Despite numerical evidence, and many theoretical progresses, the conjecture has resisted a formal proof since.
\\In the first part of this paper, we investigate a stricter version of the conjecture, expressed as follows: ''{\em Let $p_{n}$ be the n-th prime. Then, there always exist twin primes between $(p_{n}-2)^{2}$ and $p_{n}^{2}$ }''. To justify this conjecture, we formulate a prediction (based on a double-sieve method) for the number of twin prime pairs in this range, and compare the prediction with the real results for values of $p_{n}$ up to 6500000. We also analyse what should happen for higher values of $p_{n}$.
\\In the second part, we investigate the validity of the general Polignac's Conjecture. We predict the ratio of the number of solutions for any value of M divided by the number of solutions for $M = 2$, and explain how this ratio depends on the factorization of $M$. We compare the predictions with the real values for $M$ up to 3000 (and for the special case 30030) in the range of from the 1000000-th prime to the 21000000-th prime.

\end{abstract}
\part*{\large Part 1: The Twin Prime Conjecture}

\section{\normalsize Introduction}
The Twin Prime Conjecture, which is a particular case of the more general Polignac's Conjecture \cite{Polignac}) hypothesizes that there exist an infinite number of pairs $k$, $k+2$ where both $k$ and $k+2$ are prime (''twin'' primes). Despite numerical evidence (higher and higher twin prime pairs have been found), and many theoretical progresses \cite{Zhang}, \cite{Soundararajan}, the conjecture has resisted a formal proof since. 
\\
In this paper we introduce a {\bf stricter version of the conjecture}, expressed as follows: ''{\em Let $p_{n}$ be the n-th prime. Then, there always exist twin primes between $(p_{n}-2)^{2}$ and $p_{n}^{2}$ }''. We will also formulate a prediction for the number of twin prime pairs in this range.  
\\
From now on, although many arguments we present apply to the general case, we will focus on this stricter conjecture, i.e., that there exists at least one pair of twin primes between $(p_{n}-2)^{2}$ and $p_{n}^{2}$. One of the reasons to choose such a restricted range is that, as the range starts with a value higher than $p_{n-1}^{2}$,  all the prime numbers up to $p_{n}-1$ can eliminate ''candidates'' (defined in the next section) in the entire range.  This would not be true, for instance, if we used a range from 1 to $p_{n}^{2}$; in such a case, $p_{i}$ would only be significant in the portion of the range above $p_{i}^{2}$ (introducing complexity in a prediction's formulation). 
\\
We proceed along the following path:
\begin{itemize}
\item count the number of potential candidates for twin prime pairs between $(p_{n}-2)^{2}$ and $p_{n}^{2}$.
\item use a double-sieve to estimate the probability for a candidate to actually be a twin prime pair.
\item derive a prediction for the number of twin primes in the selected range.
\item compare the prediction with the real result (obtained from a computer program).
\item investigate what should happen for larger $p_{n}$ values, to show that when $p_{n}$ grows, also the predicted number of twin pairs grows.
\end{itemize}
Finally, note that for $p_{n} = 3$ and $p_{n}-2 = 1$ the solutions 3,5 and 5,7, and for $p_{n}$  = 5 and $p_{n}-2 = 3$, the solutions 11,13 and 17,19 satisfy the conjecture. In what follows we will therefore assume $p_{n} \geq 7$. 
\section{\normalsize Looking for twin primes between $(p_{n}-2)^{2}$ and $p_{n}^{2}$}
\subsection{\normalsize Counting the candidates}
Any pair of twin primes can always be expressed as $N-1$, $N+1$, where $N$ is a multiple of 6 (3 must divide $N$, otherwise it would divide one of $N-1$, $N+1$). 
\\To search for twin pairs satisfying our version of the conjecture for a given  $p_{n}$, we have to examine all the integer numbers $N_{m}= 6m$, contained in the range from $(p_{n}-2)^{2}$ to $p_{n}^{2}$. These numbers will be our ''candidates'', and for each of them we will check if both $N-1$ and $N+1$ are prime. We list these candidates in increasing order, from $N_{m\_1}$ (the smallest value above $(p_{n}-2)^{2}$  to be divisible by 6), to $N_{m\_last}$ (which is equal to $p_{n}^{2} - 7$, as explained in the next paragraph).   
\\
For all $p_{n} > 3$, $p_{n}$ modulo 3 is either 1 or 2, hence $p_{n}^{2}$ modulo 3 is 1. $N = p_{n}^{2}-1$ is an even number, and its modulo 3 is 0, hence it is divisible by 6. But we can always exclude this candidate, because N+1 is $p_{n}^{2}$, which is not a prime. Therefore the largest candidate is $p_{n}^{2}$ - 7. 
\\For a given $p_{n}$, the size of the range within which we look for twin primes is $p_{n}^{2} - (p_{n}-2)^{2}$, i.e. $4(p_{n}-1)$. If  $p_{n}-2$ is also prime, or if 3 does not divide it, $N_{m\_1}$  is equal to $(p_{n}-2)^{2}+5$, because also $(p_{n}-2)^{2}$ modulo 6 is 1. So, between $(p_{n}-2)^{2}$ and $p_{n}^{2}$ there is an ''excluded'' region of size 12 (7+5) that cannot contain useful candidates. The number of candidates will then be $(4(p_{n}-1)- 12 )/6 + 1$.\\ If $p_{n}-2$ is divisible by 3, then $(p_{n}-2)^{2}$ modulo 6 = 3; this makes $N_{m\_1}$ = $(p_{n}- 2)^{2} + 3$, and reduces the size of the excluded region to 10. In this case the number of candidates is $(4(p_{n} - 1) - 10)/6 + 1$. 
\\In both cases, when $p_{n}$ is large, a good approximation for the number of candidates is 
\\
\begin{equation}
n\_of\_candidates \cong \frac{2}{3} p_{n} \label{eq:candidates}
\end{equation}
\\
\\
\subsection{\normalsize How the values of $N$ modulo $p_{i}$ determine if $N-1$, $N+1$ is a twin prime pair}
For any candidate $N_{m}$, $N_{m} – 1$ and  $N_{m}+ 1$ are always smaller than $p_{n}^{2}$. Therefore, for any value of $m$, we know that $N_{m} – 1$ and  $N_{m}+ 1$ are prime if they are not divisible by any prime from 5 to $p_{n-1}$ (we already excluded 2 and 3, by imposing the condition that $N_{m}$ is a multiple of 6). 
\\Without loss of generality, let’s call $N$ one of these candidates. 
\\
We look at the residuals of $N$ with respects to each prime from 5 to $p_{n-1}$ . 
\begin{itemize}
\item For 5, if $N$ modulo 5 is 1, 5 divides $N-1$. If $N$ modulo 5 is 4, 5 divides $N+1$. If $N$ modulo 5 is 0, 2, or 3, $N-1$ and $N+1$ are not divisible by 5. So we discard all the $N_{m}$ values for which $N_{m}$ modulo 5 is 1 or 4 (two cases), and we keep the others (three cases). 
\item
For 7, we discard all the $N_{m}$  values for which $N_{m}$ modulo 7 is 1 or 6, (two cases) and we keep those for which modulo 7 is 0, 2,3,4, or 5 (five cases). 
\item
... 
\item
For $p_{n-1}$, we discard all the $N_{m}$ values for which $N_{m}$ modulo $p_{n-1}$ is 1 or $p_{n-1} – 1$ (two cases), and we keep all the others ($p_{n-1} – 2$ cases). 
\end{itemize}
In fact, we are applying a {\bf double-sieve} to the set of $N_{m}$ values (which form a finite arithmetic progression). 
\\
For $p_{n}$, the number of primes for which we have to check the set of $N_{m}$ values is n-3. For example, for $p_{5}$ (=
11), we have to check  $p_{4}$ (=7) and  $p_{3}$ (=5). For each of these primes $p_{j}$ , with $j$ going from 3 to $n-1$, we discard $N_{m}$ if the value of $N_{m}$ modulo $p_{j}$ is either 1 or $p_{j} - 1$. 

\subsection{\normalsize Predicting the number of twin primes}
To predict the number of prime twins for a given $p_{n}$, we assess the effect of applying our double-sieve on the set of $N_{m}$ values. In particular, we try to estimate the probability that a candidate is spared by the double sieve.  
\\
We start by evaluating the product of terms ($p_{i}-2)/p_{i}$, where the index $i$ goes from 3 ($p_{3}$ being 5) to $n-1$. We recall that, by choosing only $N_{m}$ values divisible by 6, we already took into account the effect of $p_{1}$(=2) and $p_{2}$(=3). For a given $p_{n}$, this product, which we call {\bf\em product\_pminus2\_over\_p}, consists of n-3 terms. 
\begin{equation}
product\_pminus2\_over\_p = \prod_{primes \geq p_{3}}^{p_{n-1}} \frac{p-2}{p} = \frac{3}{5}\frac{5}{7}\frac{9}{11}\frac{11}{13}\cdots\frac{p_{n-1}-2}{p_{n-1}}
\end{equation}
If we have the list of prime numbers up to $p_{n-1}$, we can get the exact value for the product; for higher values of $n$ we will have to estimate it. One can observe that the product can be rewritten as the product of $(n-2)/n$ over {\bf all} odd numbers in the range, divided by the product of $(n-2)/n$ over {\bf all composite} odd numbers, i.e.
\begin{equation}
\frac{3/5 \cdot 5/7 \cdot 7/9 \cdot 9/11 \cdot 11/13 \cdot 13/15 \cdot...\cdot (p_{n-1}-2)/p_{n-1}}{7/9 \cdot 13/15 \cdot 19/21 \cdot...\cdot (c_{n-1}-2)/c_{n-1}},
\end{equation}
where $c_{n-1}$ is the largest composite number smaller than $p_{n-1}$. In other words, we have multiplied and divided the original product by a product made with terms that correspond to all the composite numbers in the range of interest. The numerator telescopes to $3/p_{n-1}$. \\Every term in the denominator is smaller than 1; so also the denominator is smaller than 1, and it gets smaller and smaller for growing values of $n$, although less rapidly than the numerator. The denominator can be written as the product of a sequence of telescoped terms: if $p_{n}$ and $p_{n+1}$ are consecutive primes, the product coming from all the composite terms in between telescopes to $p_{n}/(p_{n+1} – 2)$. The denominator, being smaller than 1, acts as a multiplicative factor, growing with $n$, by which we multiply the $3/p_{n-1}$ value to get a first estimate of the survival probability for a candidate. The overall result is a probability value that, when $n$ grows, decreases less rapidly than $1/p_{n}$. We will then multiply this probability by the number of candidates (proportional to $p_{n}$), to find the ''expected number of twin primes'' in the interval between $(p_{n}-2)^{2}$ and $p_{n}^{2}$. 
\\
One interesting feature of this approach is to make evident that, for large $n$, the value of the denominator depends mainly on the frequency of the prime number, and only marginally on their precise distribution \footnote{As an example, we compare what would happen to a part of the product in three different hypothetical scenarios. We loop over the 4 millions odd primes from $p_{1000000}$ to $p_{5000000}$, in steps of 2. At every step we keep the right values for $p_{n}$ and $p_{n+2}$, and place an hypothetical $p_{n+1}$ in three different positions: 1), at $p_{n}+2$,  2), at $p_{n+2} -2$, and 3), at $(p_{n}+p_{n+2})/2$. We then compute $prod1$, $prod2$ and $prod3$, corresponding to the three scenarios. The first seven decimal digits of the three scenarios are identical (0.2192284), while the eight digit is 9 in the first scenario, 7 in the second, and 8 in the third scenario. So the dependency on the actual position of $p_{n+1}$ is minimal.}.
\\
If we simply multiplied the number of candidates $\cong \frac{2}{3} p_{n}$ by the numerator $3/p_{n-1}$, we would obtain a value of the order of 2, for any value of $p_{n}$. This value would be clearly insufficient to guarantee the existence of twin primes in the interval we investigate; but once we multiply this by the factor coming from the denominator, we will get a value growing with $n$. This will make the existence of twin prime couples in the selected range much more likely.  
\\According to this initial approach, the predicted number of twin primes for $p_{n}$ (i.e. between $(p_{n}-2)^{2}$ and $p_{n}^{2}$) would be 
\begin{equation}
prediction = product\_pminus2\_over\_p \cdot n\_of\_candidates \label{eq:initial}
\end{equation}
In reality, equation ~\ref{eq:initial} overestimates the number of twin primes by a factor $(1.12292)^{2}$. The reason for this derives from the following:
\begin{enumerate}
\item $\frac{p-2}{p} = \frac{p-2}{p-1} \cdot \frac{p-1}{p}$
\item The product $\prod\frac{p-1}{p}$, over all primes lower than $\sqrt{x}$, is equal to $\frac{2 \exp(-\gamma)}{log(x)}$ (Mertens' formula).
$\gamma$ is the Euler-Mascheroni constant, and the value of $2 \exp(-\gamma)$ is 1.12292. This means that the product overestimates the probability of $x$ being a prime by a factor 1.12292. 
\item The ratio between $\prod\frac{p-2}{p-1}$ and $\prod\frac{p-1}{p}$ tends to a constant value (0.66016...) called ''the twin prime constant''  when $p$ goes to infinity. The convergence is quite rapid: the value is already 0.66129 when the last term in the product is 97.
\item From 1, 2 and 3 we conclude that the product of the terms $(p-2)/p$, up to $\sqrt{x}$, overestimates the probability of $x-2$ and $x$ being a couple of twin primes by a factor $(1.12292)^{2}$.
\end{enumerate}
Therefore, our predicted value becomes
\begin{equation}
prediction = product\_pminus2\_over\_p \cdot n\_of\_candidates/(1.12292)^{2} \label{eq:prediction}
\end{equation}
\subsection{Behaviour for large numbers}
For very large numbers, we do not know exactly the number of primes and their values; therefore if we still want to make a prediction to see if the expected value of twin primes increases or decreases, in equation \ref{eq:prediction} we have to estimate the values of a) {\em product\_pminus2\_over\_p}, and of b) {\em n\_of\_candidates}. Both a) and b) require an approximation for the $n$-th prime number $p_{n}$ ; for this, we can use the first terms of an approximation \cite {Cesaro} established by Ernesto Cesàro in 1894:
\begin{equation}
p_{n} \cong n(\log n + \log \log n - 1 + (\log \log n -2)/\log n)
\end{equation}
The value of the b) term would simply be equal to 2/3 $p_{n}$, as shown earlier.
\\To estimate the a) term, we will use the formula $(log(n1)/log(n2))^2$, 
which approximates the value of the product $\prod\frac{p-2}{p}$ over all primes between $n1$ and $n2$. As this formula works better for large values, we approximate {\em product\_pminus2\_over\_p} over the primes up to {\em x} as 
%\small{
\begin{equation}
\prod_{p \leq x}\frac{p-2}{p} \cong \prod_{p \leq 10000}\frac{p-2}{p} \cdot (\log{10001}/\log{x})^2   \label{eq:approx1}
\end{equation}
%}
using the fact that we can easily compute directly the product for the primes up to 10000. The equation reduces to $\frac{C_{1}}{(\log{x})^2}$, while the number of candidates is $C_{2}x$. The final result can be written as $C\frac{x}{(\log{x})^2}$. The relevant point is that this quantity increases when x grows; this justifies the conjecture.  
\section{\normalsize Some results, to check our method and the conjecture}
For every prime number $p_{n}$ up to 250000 (corresponding to a $p_{n}^{2}$ = 62500000000) we used a computer program to test every multiple of 6 in the ranges from $(p_{n}-2)^{2}$ to $p_{n}^{2}$ against the modulo criteria expressed earlier, to find the twin prime pairs. We also examined some partial intervals (consisting of 3600000000 integers each), starting at values up to 43200000000000 (reaching values of $p_{n}$ up to more than 6500000). The highest partial interval contains only 16 squares of prime numbers.
\\
We compared the results with the predicted values derived from \ref{eq:prediction}. In these cases, as we have a complete knowledge on the primes involved, the prediction has been formulated using the computed result of the $(p-2)/p$ product. The results are shown in Figures 1 to 5.
\\
The agreement between the predicted values and the real values is quite good, and improves when $p_{n}$ gets larger. For values of $p_{n}$ above 6000000, the agreement is better than 1\%, as shown in Figure 5.
\\
\begin{figure}[hbtp]
\subfloat{
\begin{minipage}[c][1\height]{0.9\linewidth}
\centering
\includegraphics[scale=.5]{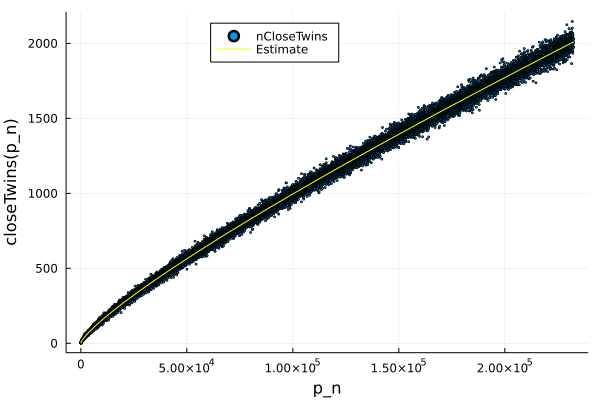}
\caption{The number of twin pairs between $(p_{n}-2)^2$ and $p_{n}^2$ for the primes up to 232367, and the prediction.}
\end{minipage} }
\hfill
\subfloat{
\begin{minipage}[c][1\height]{0.9\linewidth}
\centering
\includegraphics[scale=.5]{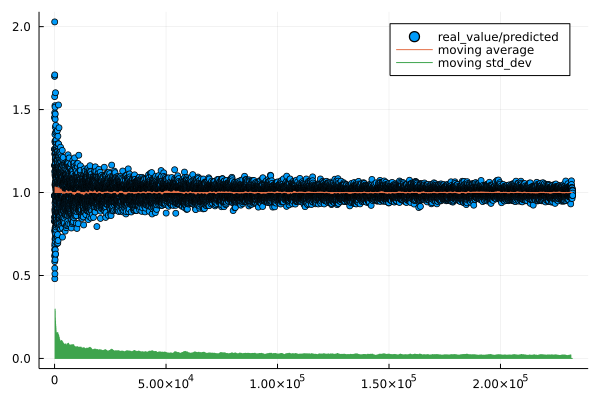}
\caption{The ratio between the real data and the prediction.}
\end{minipage} }
\end{figure}

\begin{figure}[!t]
\begin{subfigure}{.8\textwidth}
\centering
\includegraphics[scale=.5]{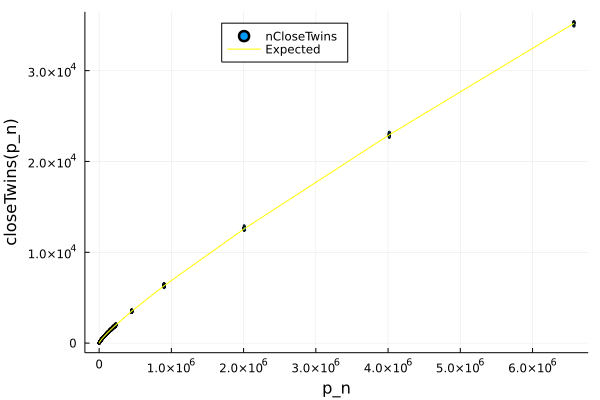}
\end{subfigure}
\caption{Same as Fig. 1 with some data for $p_{n}$ up to 6500000 , and the prediction.}
\end{figure}

\begin{figure}[!t]
\subfloat{
\begin{minipage}[c][1\width]{0.4\linewidth}
\centering
\includegraphics[scale=.25]{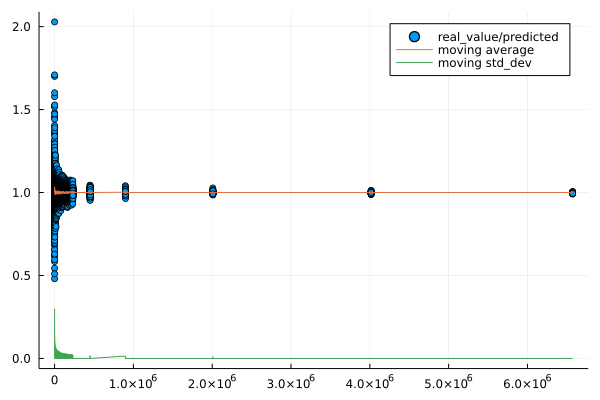}
\caption{real data/prediction}
\end{minipage} }
\hfill
\subfloat{
\begin{minipage}[c][1\width]{0.4\linewidth}
\centering
\includegraphics[scale=.25]{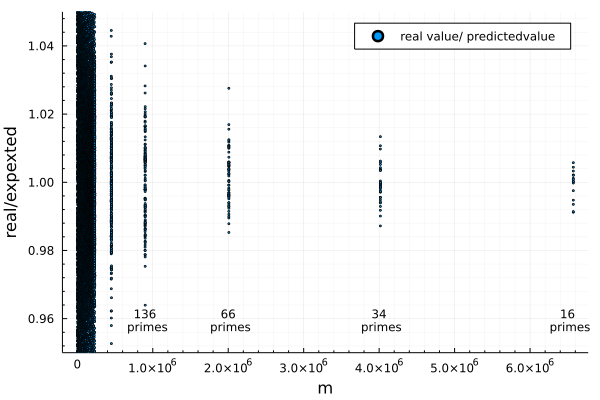}
\caption{zoomed vertically}
\end{minipage} }
\end{figure}
\section*{\normalsize Conclusion for the first part}
We have introduced a stricter version of the Twin Primes Conjecture (existence of twin pairs between $(p_{n}-2)^2$ and $p_{n}^2$) , to see if within this frame we could make reality-matching predictions. We explained in detail the process to obtain a prediction, and we have shown the good agreement between the prediction and the real data, in a range of prime values up to more than 6000000. The relative accuracy of the prediction improves with $p_{n}$; for the few prime numbers larger than 6000000, for which we counted the real number of twin pairs the agreement is better than 1\%. The expression predicting the number of twin pairs for a given $p_{n}$ increases when $n$ grows; this is at least a good indication that our version of the conjecture, and therefore also the original Twin Primes Conjecture, could be true.

\part*{\large Part 2: The Polignac's Conjecture}

The Polignac's Conjecture hypothesizes the infinite occurrence of couples of prime numbers $P$, $P + m$, where $m$ is any even number.
The second part of this article is motivated by the fact that, in many respects, the search for pair of prime numbers $P$, $P+m$ (''companion'' primes) is very similar to the search for pairs of prime numbers $P$, $P+2$ (the twin primes). We will show how, in a large enough range, a very simple relation exists between the number of occurrences of twin prime pairs and the number of ''$m$'' prime pairs. For a given large range, let's call {\em pairs\_2} the number of twin primes, and {\em pairs\_m} the number of prime couples at interval $m$. Then the ratio $pairs\_m/pairs\_2$, in first approximation, depends only on the prime numbers appearing in the factorization of $m$.
\section{\normalsize How the number of ''$m$'' prime pairs depends on the divisors of $m$}
\subsection{\normalsize Analysis}
Let's start by recalling a few facts from the first part:
\begin{itemize}
\item all the primes larger than 3 are of the form $N-1$ or $N+1$, where $N$ is a multiple of 6.
\item in a given range from $A$ to $B$ there are roughly $(B-A)/6$ multiples of 6. We can index them as $N_{i}$, with $i$ going from 1 to $max\_i$. 
\item once $B$ is defined, let's {\em max\_prime} be the largest prime number smaller than $\sqrt{B}$
\end{itemize}
Let's now start our investigation by examining different values of $m$. For a pair of integers $K$, $K+m$, we call the number $K$ the ''$Low$'' member of the pair, and $K+m$ the ''$High$'' member. We only look for pairs where $Low$ is $N_{i} \pm 1$ and $High$ is $N_{j} \pm 1$. If $m = 2$ then $i = j$; otherwise $j > i$.
\begin{itemize}
\item $m$ = 2. In this case, as we have seen earlier, the $Low$ member is $N_{i}-1$, and the $High$ member is $N_{i}+1$. Within the defined range, we have $max\_i$ possible pairs. If we choose a single value for $i$ such that $N_{i}-1$ is prime, the probability that also $N_{i}+1$ is prime is proportional\footnote{the product must be divided by 1.12292} to the product $\prod (p-2)/(p-1)$, over {\bf all} primes $p$ from 5 to $max\_prime$. Let's call this single probability ''$sp(2)$''. The expected number of twin pairs in the range is $max\_i \cdot sp(2)$.

\item $m$ = 4, or a power of 2. Take $m$ = 4. Now the low member is $N_{i}+1$, and the high member  $N_{i+1}-1$. Both the number of possible pairs and the single probability are the same as for $m$ = 2. Therefore we expect a number of {\em pairs\_4} roughly equal to the number of twin pairs.
\\For $m$ = 8, 32, 128... the situation is the same as for $m$ = 2; the only difference is that in the formula for the $High$ member $N_{i+k}+1$, $k$ takes larger and larger values\footnote{This could have a negative impact on the probability if, due to the higher values of $High$, new primes entered in the computation of the product for its single probability. We assume that a good choice of the range's size and starting point minimizes this effect. We also keep the upper limit of the range open for the $High$ values, so to still have $max\_i$ candidates.}. \\Likewise, for $m$ = 16, 64, 256... the situation is the same as for $m$ = 4. 
\\For all these cases $sp(m)$ = $sp(2)$, and the number of candidates is the same. The expected number of $pairs\_m$ is the same as the expected number of twin pairs.
\item $m$ = 6. Here, for every $N_{i}$, we have {\bf two} possibilities: $Low$ = $N_{i}-1$, $High$ = $N_{i+1}-1$, {\bf and} $Low$ = $N_{i}+1$, $High$ = $N_{i+1}+1$. The single probability is still the same as for $m$ = 2. But, because we have twice the number of candidates, the expected number of {\em pairs\_6} is equal to {\bf two times} the number of twin pairs. The same is true for any value of $m$ whose only divisors are 2 and 3 (12, 18, 24, 36, 48...).
\item $m$ = 10. Here $Low$ is $N_{i}+1$, and $High$ is $N_{i+2}-1$. For each value of $i$, we have one candidate. The notable difference with the previous cases is that if $Low$ is prime, $High$ is not divisible by 5. Therefore we have to remove the term (5-2)/(5-1) = 3/4 from the product for the single probability of $High$ to be prime. So the full $(p-2)/(p-1)$ is multiplied by 4/3 to get the single probability $sp(10)$, and the expected number of {\em pairs\_10} is 4/3 the number of twin pairs. This is a direct consequence of the fact that $m$ is a multiple of 5. This applies also to $m$ = 20, 40, 50, 80, 100, ...
\item $m$ = 14. This is a multiple of 7. It is almost exactly like the case $m$ = 10; the only difference is that the full product is multiplied by 6/5 (and not 4/3). The expected number of {\em pairs\_14} is 6/5 the number of twin pairs.
\item $m$ = 30 (= 5 x 6). Here we combine the features of $m$ = 6 and of $m$ = 5. Compared to $m$ = 2, we have twice the number of candidates (from $m$ = 6), and the single probability of $m$ = 10, i.e. 4/3 larger than the one for $m$ = 2. Multiplying these factors, the expected number of {\em pairs\_30} is 8/3 the number of twin pairs. The same is true for $m$ = 60, 90, 120, 150, 180.
\item $m$ = 210 (= 5 x 6 x 7). Compared to $m$ = 30, we now also have 7 in the factorization of $m$. The single probability for $m$ = 30 must be multiplied by 6/5, which brings the expected number of {\em pairs\_210} to 16/5 the number of twin pairs.
\item and so on...
\end{itemize}
To summarize: once the range in which we are searching for companion primes is defined
\begin{itemize}
\item if $m$ is 2, or a power of 2, the single probability $sp(2)$, i.e. the probability that if $N$ is prime also $N+2$ is prime is proportional to the product $\prod (p-2)/(p-1)$, over {\bf all} primes $p$ from 5 to $max\_prime$.
\\We call $sp(2)$ the ''base single probability'', because, as we have shown, $sp(2)$ is the smallest of all, and all the other $sp(m)$ values are related to it. \\Once the range to search is defined, $sp(2)$ depends only on the prime numbers whose squared values are smaller than the high limit of the range. Ideally, we could select the range in such a way it does not contain any $p_{n}^{2}$ value, so that $sp(2)$ is constant over the entire range. But even if this was not the case, as we are mainly interested in the relations between $sp(2)$ and the other $sp(m)$ values, this would not be a problem, because they would be affected in the same way. Also, a few very large $p_(n)$ values do not affect significantly the value of the product.
\item if $m$ is a multiple of 6 (i.e. if 3 divides $m$), {\bf the number of candidates doubles}. As $sp(m)$ does not depend on the fact that $m$ is divisible by 3, the total effect is to multiply by a factor two the number of expected pairs.
\item For every other divisor $p$ of $m$, with $p \geq 5$, the single probability $sp(2)$ must be multiplied by $(p-1)/(p-2)$ to obtain $sp(m)$.
\end{itemize}
\subsection{\normalsize Numerical results}
We examined a range of 20000000 primes, from 15485867 (the million-th odd prime) to 393342743 (the 21-million-th odd prime). In this range, we  counted all the pairs of primes $N$, $N+m$, for all even values of $m$ from 2 to 3000, and for the ''special'' case 30030 (2 x 3 x 5 x 7 x 11 x 13). There are 1399293 {\em pairs\_2} (twin pairs), 1400189 {\em pairs\_4}, 2798981 {\em pairs\_6}, \dots, 5429051 {\em pairs\_30030}.
\\For every $m$ we divide {\em pairs\_m} by {\em pairs\_2}, to obtain the ''occurrence ratio'', and compute the ''expected ratio'', according to the criteria explained above. In our sample of $m$-values, the expected ratio goes from 1.0 (for $m$ = 2 and all powers of 2) to 3.8787 (for $m$ = 30030).
\\For every $m$ we divide the occurrence ratio by the expected ratio. For every $m$ we obtain results very close to 1 (the minimum is 0.9989, and the maximum 1.0020. The average is 1.0005, and the standard deviation is 0.0004). This shows the correctness of the model.
Figures 6 and 7 show these results (we have excluded the point $m$ = 30030, to avoid a large hole in the plot).\\Figure 6 shows the number of $pairs\_m$. It is clearly possible to distinguish the base line, a line for those $m$ values divisible by 3, a line for those divisible by 3 and 5, one for those divisible 3, 5 and 7, and the highest point, corresponding to $m$ = 2310 = 2x3x5x7x11. More on the right, we find the second highest point, for $m$ = 2730, = 2x3x5x7x13. The base line looks more crowded than expected; this is due to all the $m$ values whose factorization contains only 2 and a much larger number, which has a very limited impact on the number of prime pairs. The same applies to the "divisible by 3" line, and to a lesser extent to the "divisible 3 and 5" line.
Figure 7 shows how the numbers of $pairs\_m$ found are in very good agreement with the predictions.
\begin{figure}[!h]
\subfloat{
\begin{minipage}[c][1\height]{0.9\linewidth}
\centering
\includegraphics[scale=.4]{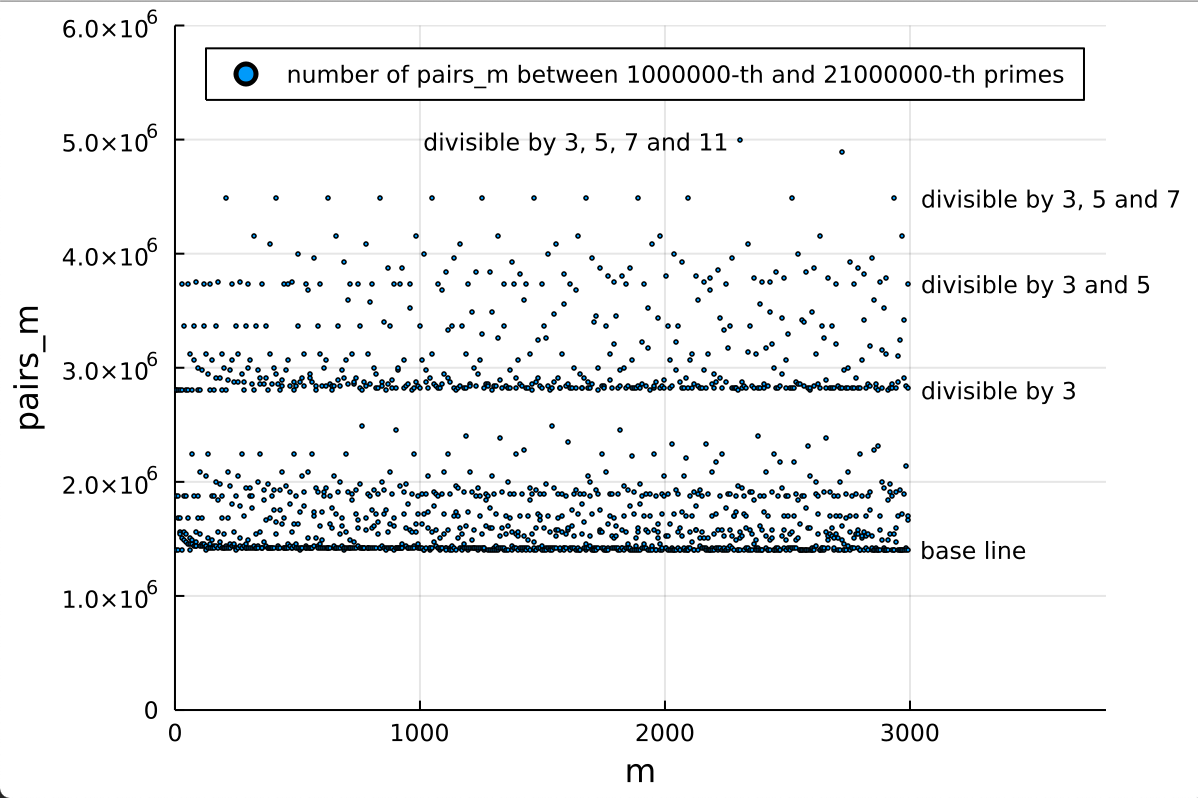}
\caption{The number of pairs of primes $N$, $N+m$}
\end{minipage} }
\end{figure}
\hfill
\begin{figure}[!h]
\subfloat{
\begin{minipage}[c][1\height]{0.9\linewidth}
\centering
\includegraphics[scale=.4]{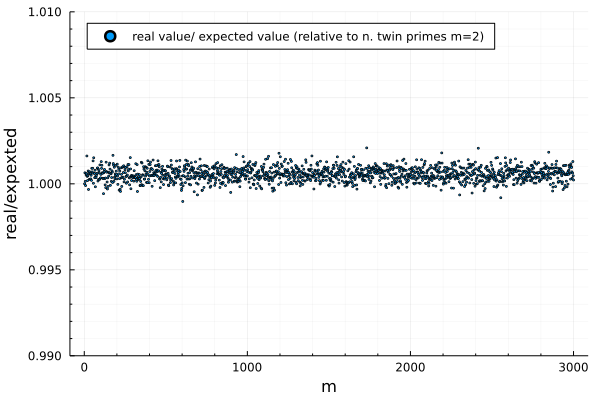}
\caption{Comparison between real and expected number of prime pairs occurrences}
\end{minipage} }
\end{figure}
\\

\section*{\normalsize Conclusion}
We have shown the existence of a relation between the frequencies of prime pairs at intervals $m_{1}$ and $m_{2}$, and how to quantify the ratio between these frequencies. Based on this, we can make another consideration. In this process, we predicted and compared the frequencies for every value of $m$ against the frequency for $m$ = 2, which we used as reference. But there is nothing special about the value $m$ = 2. We could have chosen any value of $m$, for example $m_{x}$ as reference, and then predicted the occurrence frequencies for all the other values of m relative to the occurrence frequency for $m_{x}$. In \cite{Zhang}, Y. Zhang has proved the existence of (at least) one unknown value $m_{unk}$, between 2 and 70000000, such that there exist an infinity number of prime pairs at interval $m_{unk}$. Further work by J. Maynard \cite{Maynard_1} has first lowered the upper limit for $m_{unk}$ to 600, and then, in collaboration with the polymath project, to 246 \cite{Polymath}. We could have selected this $m_{unk}$ as reference. Then, using the method discussed in the last section, we would obtain occurrence frequencies larger than zero for any $m$ value; this would ''prove'' the Polignac's Conjecture, and, as a special case, the Twin Primes Conjecture.

\bibliography{Primes}
\bibliographystyle{unsrt} 
\end{document}